\documentclass{amsart}

\usepackage[psamsfonts]{amssymb}

\theoremstyle{plain}
\newtheorem{theorem}{Theorem}[section]

\newtheorem{cor}[theorem]{Corollary}

\theoremstyle{definition}
\newtheorem{defn}{Definition}[section]

\theoremstyle{remark}

\title[Semiclassical multivariate
 orthogonal polynomials]{A semiclassical perspective on multivariate
 orthogonal polynomials}

\author[M. \'{A}lvarez de Morales]{Mar\'{\i}a \'{A}lvarez de Morales}
\address[M. \'{A}lvarez de Morales]{Departamento de Ma\-te\-m\'a\-ti\-ca Apli\-ca\-da,
Uni\-ver\-si\-dad de Gra\-na\-da,
 Gra\-na\-da,
Spain} \email{alvarezd@ugr.es}

\author[L. Fern\'andez]{Lidia Fern\'andez}
\address[L. Fern\'andez]{Departamento de Ma\-te\-m\'a\-ti\-ca Apli\-ca\-da,
Uni\-ver\-si\-dad de Gra\-na\-da,
 Gra\-na\-da,
Spain} \email{lidiafr@ugr.es}

\author[T. E. P\'erez]{Teresa E. P\'erez}
\address[T. E. P\'erez]{Departamento de Matem\'atica Aplicada,
  and Instituto Carlos I de F\'{\i}sica Te\'orica y Computacional,
  Universidad de Granada,
Granada,
 Spain}
\email{tperez@ugr.es}

\author[M. A. Pi\~{n}ar]{Miguel ~~ A. ~~ Pi\~{n}ar}
\address[M. A. Pi\~{n}ar]{Departamento de Matem\'atica Aplicada,
  and Instituto Carlos I de F\'{\i}sica Te\'orica y Computacional,
   Universidad de Granada,
Granada,
 Spain}
\email{mpinar@ugr.es}

\thanks{Partially supported by Ministerio de Ciencia y Tecnolog\'{\i}a
        (MCYT) of Spain and by the European Regional Development Fund
        (ERDF) through the grant MTM 2005--08648--C02--02, and Junta de
        Andaluc\'{\i}a, Grupo de Investigaci\'on FQM 0229.}

\begin{document}



\begin{abstract}
Differential properties for orthogonal polynomials in several
variables are studied. We consider multivariate orthogonal
polynomials whose gradients satisfy some quasi--orthogonality
conditions. We obtain several characterizations for these
polynomials including  the analogous of the semiclassical Pearson
differential equation, the structure relation and a
differential--difference equation.
\end{abstract}

\maketitle

\section{Introduction}

In the univariate case, a \emph{semiclassical moment functional} $u$
is a quasi--definite functional satisfying a distributional Pearson
equation
\begin{equation}\label{peq}
     D(\phi u)= \psi u
\end{equation}
where $\phi$ and $\psi$ are polynomials with $\deg(\psi) \ge 1$.
They constitute a natural extension of classical functionals
(Hermite, Laguerre, Jacobi, and Bessel) and they have been
extensively analyzed during the last two decades (see \cite{HR},
\cite{Ma}). Obviously, the sequence $\{p_n\}_n$ of orthogonal
polynomials associated with $u$ is also called \emph{semiclassical}.

The distributional Pearson equation (\ref{peq}) plays a key role in
the study of differential properties of semiclassical orthogonal
polynomials. In fact, this equation allows us to characterize
semiclassical polynomials as the only sequences of orthogonal
polynomials satisfying one of the following equivalent properties:

\begin{enumerate}
\item[(a)] the so--called \emph{structure relation},

\item[(b)] the quasi--orthogonality of the derivatives,

\item[(c)] the second order differential--difference relation.

\end{enumerate}
For further properties and characterizations see, for instance,
\cite{Ma}.

In the multivariate case, we will call \emph{semiclassical} a
quasi--definite moment functional $u$ satisfying the matrix {\it
Pearson--type} equation
$$\hbox{\rm div~} (\Phi ~u) =
\Psi^t ~u,
$$
where $\Phi$ is a $d \times d$ symmetric polynomial matrix and
$\Psi$ is a $d \times 1$ polynomial matrix with $\deg \Psi \ge 1$,
and such that $ \det\langle u, \Phi\rangle \neq 0.$ Of course, this
definition includes all the ``classical" moment functionals in the
usual literature (\cite{AK,KKL1,KKL2,Koo,KS,Li,Su}).

The aim of our contribution is to analyze the extension of
characterizations (a), (b), and (c) to the multivariate case by
means of a matrix formalism. In the bivariate case, the first of
these three characterizations was the main result of a previous
publication (see \cite{AFPP1}). The two other characterizations
constitute the main objective of this work.

The structure of the paper is the following. In Section 2, a basic
background about moment functionals and multivariate orthogonal
polynomials is given in order to allow the reader to be familiar
with such concepts. In Section 3 we present two new
characterizations for sequences of multivariate semiclassical
orthogonal polynomials. Theorems \ref{teor3} and \ref{teor4} are,
respectively, the extensions of the quasi--orthogonality of the
derivatives and the differential--difference relation satisfied by
semiclassical orthogonal polynomials in one variable. Finally, in
Section 4 we explore two non--trivial examples of multivariate
semiclassical orthogonal polynomials.

\section{Background}

We denote by ${\mathcal P}$ the linear space of polynomials in $d$
variables with real coefficients, and by ${\mathcal P}'$ its
topological dual (see \cite{Ma}). Let ${\mathcal P}_n$ be the set
of real polynomials of total degree not greater than $n$, and
$\Pi_n$ the space of homogeneous polynomials of degree $n$ in $d$
variables. It can be observed that
$$
\dim {\mathcal P}_n = \binom{n+d}{n} \qquad \mbox{and}\qquad
\dim{\Pi}_n=\binom{n+d-1}{n} = r_n.
$$
Let ${\mathcal M}_{h\times k}(\mathbb{R})$ be the linear space of
$(h\times k)$ real matrices, and ${\mathcal M}_{h\times
k}({\mathcal P})$ the linear space of $(h\times k)$ polynomial
matrices. In addition, $I_h$ represent the identity matrix of
order $h$.

If  $M= \left(m_{i,j}(x,y)\right)_{i,j=1}^{h,k}$ represents a
$(h\times k)$ polynomial matrix, we define the {\it degree of} $M$
by
$$\deg M = \max\{\deg m_{i,j}(x,y), 1\le i\le h, 1\le j\le
k\}\ge 0.$$

Now, we review some basic definitions and tools about multivariate
orthogonal polynomials that we will need in the rest of the paper.
For a complete description of this and other related subjects see,
for instance, \cite{AFPP1,DX,FPP3,Ko1,Ko2,Su,Xu1}.

Let $\mathbb{N}_0$ be the set of nonnegative integers. For $\alpha
= (\alpha_1,\dots,\alpha_d)\in\mathbb{N}_0^d$ and
$\textrm{x}=(x_1,\dots,x_d)\in\mathbb{R}^d$ we write
$\textrm{x}^{\alpha}=x_1^{\alpha_1}\cdots x_d^{\alpha_d}$. The
number $|\alpha| = \alpha_1 + \dots + \alpha_d$ is called the
\emph{total degree} of $\textrm{x}^{\alpha}$.

Let $\{\mu_{\alpha}\}_{\alpha\in \mathbb{N}_0^d}$ be a
multi--sequence of real numbers and let us denote by $u$, the only
linear functional in ${\mathcal P'}$ satisfying
$$
\langle u,\textrm{x}^{\alpha}\rangle =\mu_{\alpha}.
$$
where, as usual, $\langle \, , \, \rangle$ stands for the duality
bracket. Then, $u$ is called the \emph{moment functional} determined
by $\{\mu_{\alpha}\}_{\alpha\in \mathbb{N}_0^d}$, and the number
$\mu_{\alpha}$ is called the \emph{moment of order} $\alpha$.

Let $(P_{\alpha}^n)_{|\alpha|=n}$ be a sequence of polynomials of
total degree $n$. Using the matrix notation introduced in
\cite{Ko1,Ko2} and \cite{Xu1}, we can denote by $\mathbb{P}_n$ the
vector polynomial
$$\mathbb{P}_n=(P_{\alpha}^n)_{|\alpha|=n}=
(P_{\alpha^{(1)}}^n,P_{\alpha^{(2)}}^n,\cdots,P_{\alpha^{(r_n)}}^n)^t
\in {\mathcal M}_{{r_n}\times 1}({\mathcal P}_n),$$ where
$\alpha^{(1)},\alpha^{(2)},\cdots,\alpha^{(r_n)}$ are the elements
in $\{\alpha \in \mathbb{N}_0^d:\,|\alpha|=n \}$ arranged
according to the lexicographical order.

 When $\{\mathbb{P}_m\}_{m=0}^n$
is a basis of ${\mathcal P}_n$ for each $n\geq 0$, then
$\{\mathbb{P}_n\}_{n\ge0}$ is called a {\it polynomial system}
(PS).

\begin{defn}
We will say that a PS  $\{\mathbb{P}_n\}_{n\geq 0}$ is a
\emph{weak orthogonal polynomial system (WOPS)} with respect to
$u$ if
\begin{equation}\label{oc}
    \begin{array}{l}
   \langle u,{\mathbb{P}}_{n}{\mathbb{P}}_{m}^t\rangle =0, \quad m \neq n, \\
   \langle u,{\mathbb{P}}_{n}{\mathbb{P}}_{n}^t\rangle =H_n,
   \end{array}
\end{equation}
where $H_n\in {\mathcal M}_{r_n}(\mathbb{R})$ is a non singular
matrix.
\end{defn}

This definition means that every polynomial component in
$\mathbb{P}_n$, is orthogonal to all the polynomials of lower
degree, but two polynomial components of the same degree don't
have to be orthogonal.

A moment functional $u$ is said \emph{quasi--definite} if there
exists a WOPS with respect to $u$ (\cite{KKL1}).

Now, we are going to recover the three term recurrence relation for
orthogonal polynomials in several variables (\cite{DX,Ko1,Ko2,Xu1}).
This relation takes a vector-matrix form and it plays an essential
role in understanding the structure of orthogonal polynomials, as in
the univariate case.

For $n\geq 0$, there exist unique matrices $A_{n,i}\in {\mathcal
M}_{r_n\times r_{n+1}}(\mathbb{R})$, $B_{n,i}\in {\mathcal
M}_{r_n\times r_n}(\mathbb{R})$ and $C_{n,i}\in {\mathcal
M}_{r_n\times r_{n-1}}(\mathbb{R})$, $i=1,2, \ldots, d$, such that

\begin{equation}\label{rr3t}
   x_i\,\mathbb{P}_n= A_{n,i}\,\mathbb{P}_{n+1}+B_{n,i}\,\mathbb{P}_n
   + C_{n,i}\,\mathbb{P}_{n-1}, \quad 1\le i\le d,
\end{equation}
where $\mathbb{P}_{-1}= 0$ and $C_{-1,i}=0$.

Moreover, $\mbox{rank} \, A_{n,i} = \mbox{rank}\, C_{n+1,i}^t =
r_n $ and $
    \mbox{rank}\, A_n = \mbox{rank}\, C_{n+1}^t = r_{n+1},
$
where $A_n=\left(\begin{array}{c}
  A_{n,1} \\
  \vdots \\
  A_{n,d}
\end{array}\right )$ and $C_{n+1}^t=\left(\begin{array}{c}
  C_{n+1,1}^t \\
  \vdots\\
  C_{n+1,d}^t \\
\end{array}\right ).$

As a consequence of the three term recurrence relation, it is
possible to obtain the polynomial $\mathbb{P}_{n+1}$ in terms of the
vector polynomials $\mathbb{P}_n$ and $\mathbb{P}_{n-1}$ (see
\cite{DX}).

Let $D_n^t=(D_{n,1}^t \,\ldots \, D_{n,d}^t )\in {\mathcal
M}_{r_{n+1}\times d r_n}(\mathbb{R})$ be a left generalized inverse
of $A_n$. There exist matrices $E_n^{n+1}$ and $E_{n-1}^{n+1}$ such
that
\begin{equation}\label{rr3t1}
\mathbb{P}_{n+1}=\sum_{i=1}^d x_i\, D_{n,i}^t \,\mathbb{P}_n
 + E_n^{n+1} \,\mathbb{P}_n +E_{n-1}^{n+1} \, \mathbb{P}_{n-1}.
 \end{equation}

Moreover, since the matrices $C_{n,i}$, for $i=1,2, \ldots, d,$ are
full rank, we can obtain $\mathbb{P}_{n-1}$ in terms of the vector
polynomials $\mathbb{P}_{n+1}$ and $\mathbb{P}_n$. In fact, using
the left generalized inverse of $C_{n,i}$, that we will denote by
$G_{n,i}$, we deduce the following relations:

\begin{equation}\label{rr3t2}
\mathbb{P}_{n-1}=- G_{n,i}A_{n,i} \,\mathbb{P}_{n+1}+ (x_i\,
G_{n,i}
 - G_{n,i} B_{n,i})\, \mathbb{P}_{n}, \quad 1\le i\le d.
\end{equation}

\section{Semiclassical multivariate orthogonal polynomials}

First, we introduce the concept of semiclassical moment functional.
For $d=2$, this definition was given in \cite{AFPP1}.

\begin{defn} \label{d1}
A quasi--definite moment functional $u$ is said to be
\emph{semiclassical} if it satisfies the matrix {\it
Pearson--type} equation
\begin{equation}\label{tpm}
\hbox{\rm div~} (\Phi ~u) = \Psi^t ~u,
\end{equation}
where
$$
    \Phi = \left(\phi_{ij}\right)_{i,j=1}^d \in {\mathcal M}_{d}({\mathcal
                P}),\qquad
    \Psi =\left(\psi_i\right)_{i=1}^d \in {\mathcal M}_{d\times 1}({\mathcal
            P}),
$$
are polynomial matrices such that $\Phi$ is symmetric, $\deg \Phi
= p\ge 0$, $\deg \Psi = q \ge1$, and
\begin{equation}\label{phi1}
 \det\langle u, \Phi\rangle \neq
0. \end{equation} We denote by $s=\max\{p-2, q-1\}\ge 0.$
\end{defn}

Expression (\ref{tpm}) means
$$\langle \hbox{\rm div~} (\Phi ~u), f\rangle  = \langle \Psi^t
~u,f\rangle,$$ that is, $$\langle u, \Phi ~ \nabla f + \Psi ~
f\rangle = 0, \quad\forall f\in \mathcal{P}.$$

The natural extension of the above property for matrices involves
the Kronecker product (see, for instance \cite{HJ}, p. 242),
\begin{equation}\label{tpm2}
\hbox{\rm div~} ((\Phi\otimes I_h) ~u) = (\Psi^t\otimes I_h) ~u,
\qquad h\ge 1.
\end{equation}
Relation (\ref{tpm2}) is equivalent to
\begin{equation}\label{tilde}
(\Phi\otimes I_h)~\nabla~u = (\tilde\Psi\otimes I_h) ~u, \qquad
h\ge1,
\end{equation}
where $\tilde\Psi = \Psi - (\hbox{\rm div~}\Phi)^t$.

\noindent {\bf Remark} If $s=0$, that is, $\deg\Phi=p\le 2$ and
$\deg\Psi =1$, we recover the definition of classical WOPS given
in \cite{FPP3}, which includes the classical bivariate orthogonal
polynomials studied by H. L. Krall and I. M. Sheffer (\cite{KS}),
and other authors (\cite{DX,KKL1,KKL2,Li,Su}).

A WOPS with respect to a semiclassical moment functional $u$ is
called \emph{semiclassical}.

Now, we are going to prove three characterizations for multivariate semiclassical
orthogonal polynomials: structure relation, quasi--orthogonality
relation for the gradients, and a differential--difference relation.
From now on, we will denote by $\{\mathbb{P}_n\}_{n\geq 0}$ a given
WOPS associated with a quasi--definite moment functional $u$.

\begin{theorem}[Structure relation]\label{estructura} The moment functional $u$ is
semiclassical if and only if $\{\mathbb{P}_n\}_{n\geq 0}$
satisfy
\begin{equation}\label{sr}
        \Phi~\nabla \mathbb{P}_n^t = \sum_{j=n-s-1}^{n+p-1} (I_d\otimes \mathbb{P}_j^t) F_j^n,
        \quad {\hbox {\rm for}} \quad n\ge s+1,
    \end{equation}
 where $F_{j}^n\in {\mathcal M}_{d\,r_j\times r_n}(\mathbb{R})$.

\end{theorem}

\begin{proof} See Theorem 5 in \cite{AFPP1}.
\end{proof}

Using relations (\ref{rr3t1}) and (\ref{rr3t2}), a shorter structure
relation whose coefficients are polynomial matrices can be deduced.

\begin{cor}
If $u$ is semiclassical, then $\{\mathbb{P}_n\}_{n\geq 0}$ satisfy
    $$
   \Phi~\nabla \mathbb{P}_n^t = (I_d\otimes
   \mathbb{P}_{n+1}^t) M_{1}^n + (I_d\otimes
  \mathbb{P}_{n}^t) M_{2}^n,\quad {\hbox {\rm for}} \quad n\ge s+1,
    $$
where $M_{i}^n$ are polynomial matrices with $\deg(M_{1}^n)\leq s$
and $\deg(M_{2}^n)\leq s + 1$.
\end{cor}

From Theorem \ref{estructura}, we can obtain the second
characterization for semiclassical orthogonal polynomials.

\begin{theorem}[Quasi--orthogonality relation for gradients]\label{teor3}

The func\-tio\-nal $u$ is semiclassical if and only
if, for $n\ge s+1$,
$$
\langle u,  (\nabla \mathbb{P}_m^t)^t \, \Phi~\nabla \mathbb{P}_n^t
\rangle = 0,
        \qquad \qquad 0\le m < n-s.
$$
\end{theorem}

\begin{proof} Suppose that $u$ is semiclassical.
From (\ref{sr}), we can write
\begin{eqnarray*}
    \langle u, (\nabla \mathbb{P}_m^t)^t \Phi \nabla \mathbb{P}_n^t\rangle
    = \sum_{j=n-s-1}^{n+p-1}
\langle
    u, (\nabla \mathbb{P}_m^t)^t (I_d\otimes
\mathbb{P}_j^t)\rangle \, F_{j}^n.
\end{eqnarray*}
Taking into account that
$$
(\nabla \mathbb{P}_m^t)^t(I_d\otimes \mathbb{P}_{j}^t)=\left(
\begin{array}{c|c|c|c}
  (\partial_1\mathbb{P}_m)\mathbb{P}_j^t & (\partial_2\mathbb{P}_m)\mathbb{P}_j^t
  & \cdots & (\partial_d\mathbb{P}_m)\mathbb{P}_j^t \\
\end{array}
\right),
$$
we obtain
$$
\langle u, (\nabla \mathbb{P}_m^t)^t(I_d\otimes
\mathbb{P}_{j}^t)\rangle = 0, \qquad m - 1 < j,
$$
that is, $$ \langle u, (\nabla \mathbb{P}_m^t)^t \Phi \nabla
\mathbb{P}_n^t\rangle = 0, \qquad m < n -s.$$

\bigskip

Reciprocally, assume that the quasi--orthogonality relations hold,
we define
$$\Psi=-\displaystyle{\sum_{i=0}^{s+1}}
\langle u, (\nabla\mathbb{P}_1^t)^t \, \Phi\,
\nabla\mathbb{P}_i^t\rangle H_i^{-1}\, \mathbb{P}_i.$$

Since $\nabla\mathbb{P}_1^t=I_2$, we obtain
$$\begin{array}{ll}
    \langle \hbox{\rm div~}(\Phi~u),\mathbb{P}_n^t\rangle  =&
     - \langle u,\Phi\,\nabla\mathbb{P}_n^t\rangle =
    -\langle u,(\nabla\mathbb{P}_1^t)^t\,\Phi\,\nabla\mathbb{P}_n^t\rangle.
\end{array}
$$
On the other hand,
$$
    \langle \Psi^t~u,\mathbb{P}_n^t\rangle  =
     \langle u,\Psi\,\mathbb{P}_n^t\rangle =
    - \sum_{i=0}^{s+1}\langle u,
(\nabla\mathbb{P}_1^t)^t\, \Phi\, \nabla\mathbb{P}_i^t\rangle\,
H_i^{-1} \langle u,\mathbb{P}_i \, \mathbb{P}_n^t\rangle.
$$

If $n\ge s+2$, using the above relations, we get:
$$\langle \hbox{\rm div~}(\Phi~u),\mathbb{P}_n^t\rangle  =0 = \langle \Psi^t~u,\mathbb{P}_n^t\rangle.
$$
Furthermore, for $0\le n\le s+1$,
$$
\langle \hbox{\rm div~}(\Phi~u),\mathbb{P}_n^t\rangle=-\langle
u,(\nabla\mathbb{P}_1^t)^t\,\Phi\,\nabla\mathbb{P}_n^t\rangle =
\langle \Psi^t~u,\mathbb{P}_n^t\rangle,
$$
and then, the result follows.
\end{proof}

Now, we deduce a matrix differential--difference relation for
semiclassical orthogonal polynomials in several variables.

Let us define the differential operator
$$
L[f]=\hbox{\rm div~}(\Phi \nabla f) + \tilde{\Psi}^t \nabla f,
\qquad \forall \, f \in{\mathcal P}.
$$
 Therefore,
the Lagrange adjoint of $L$ is given by
\begin{equation}\label{Ladjunto}
L^*[u]=\hbox{\rm div~}(\Phi \nabla u)- \hbox{\rm
div~}(\tilde{\Psi} u),
\end{equation}
since it satisfies
$$\langle L^*[u], f\rangle = \langle u, L[f] \rangle, \quad \forall \, f \in{\mathcal P}.$$

Using the explicit expression of the polynomial matrices $\Phi$
and $\Psi$, the operator $L$ can be written as follows

$$L[f]= \sum_{i,j=1}^d \phi_{ij}\,\partial^2_{ij}f + \sum_{i=1}^d
\psi_i\,
\partial_i f,
$$
and, from the above expression, we deduce $\deg L[f]\le s+\deg f
$.

\begin{theorem}[Matrix differential--difference relation] \label{teor4} A
func\-tional $u$ is semiclassical if and only if there exist matrices
$\Lambda_i^n\in {\mathcal M}_{r_i\times r_n}(\mathbb{R})$, such that
\begin{equation}\label{pip}
L[\mathbb{P}_n^t] =\sum_{i=n-s}^{n+s}\mathbb{P}_i^t \Lambda_i^n,
\qquad n\ge s+1.
\end{equation}
When $n\le s$, relation (\ref{pip}) reads
$$L[\mathbb{P}_n^t] =\sum_{i=1}^{n+s}\mathbb{P}_i^t \Lambda_i^n,$$
that is, $\Lambda_0^n=0, \forall n\ge 0.$

\end{theorem}

\begin{proof} If $u$ is a semiclassical functional, Lemma 4.1 in \cite{FPP3} provides
$$
\langle u, \mathbb{P}_m~\hbox{\rm div~}(\Phi \nabla
\mathbb{P}_n^t)\rangle = \langle u,  \hbox{\rm div~}((I_d\otimes
\mathbb{P}_m)~ \Phi~\nabla \mathbb{P}_n^t) \rangle - \langle u,
(\nabla \mathbb{P}_m^t)^t ~ \Phi~ \nabla \mathbb{P}_n^t\rangle,$$
for $n, m\ge0$. Besides, from (\ref{tilde}), we deduce
\begin{eqnarray*}
\langle u,  \hbox{\rm div~}((I_d\otimes \mathbb{P}_m) \Phi\nabla
\mathbb{P}_n^t) \rangle &=&-\langle \nabla u, (\Phi\otimes
I_{r_m})(I_d\otimes \mathbb{P}_m)
 \nabla \mathbb{P}_n^t \rangle \\
& = & - \langle (\Phi\otimes I_{r_m}) \nabla u,  (I_d\otimes
\mathbb{P}_m)
 \nabla \mathbb{P}_n^t \rangle\\
& = & -\langle (\tilde{\Psi}\otimes I_{r_m}) u, (I_d\otimes
\mathbb{P}_m)
 \nabla \mathbb{P}_n^t \rangle\\
 &=& -\langle  u,
(\tilde{\Psi}^t\otimes I_{r_m})(I_d\otimes \mathbb{P}_m)
 \nabla \mathbb{P}_n^t \rangle\\
& = & -\langle  u,  \mathbb{P}_m\tilde{\Psi}^t
 \nabla \mathbb{P}_n^t \rangle.
\end{eqnarray*}
Therefore, we have
\begin{equation}\label{12b}
\langle u, \mathbb{P}_m L[\mathbb{P}_n^t]\rangle = - \langle u,
(\nabla \mathbb{P}_m^t)^t \Phi\, \nabla \mathbb{P}_n^t\rangle.
\end{equation}

Now, since $L[\mathbb{P}_n^t]$ is a polynomial matrix of degree at
most $n+s$, we obtain the expansion
$$L[\mathbb{P}_n^t]=
\sum_{i=0}^{n+s}\mathbb{P}_i^t \Lambda_i^n, $$ where $\Lambda_m^n$
is given by
$$\langle u, \mathbb{P}_m L[\mathbb{P}_n^t]\rangle= \langle u,
\mathbb{P}_m \sum_{i=0}^{n+s}\mathbb{P}_i^t \Lambda_i^n \rangle= H_m
\Lambda_m^n,$$
 and using (\ref{12b}), we get
$$H_m
~\Lambda_m^n =  -
\langle u, (\nabla \mathbb{P}_m^t)^t~ \Phi ~ \nabla
\mathbb{P}_n^t\rangle.$$
In the case $n\ge s+1$, from Theorem \ref{teor3},
$$H_m ~\Lambda_m^n =0, \quad 0\le m< n-s,$$
so, relation (\ref{pip}) follows.

When $n\le s$,
$$H_0 ~ \Lambda_0^n = -
\langle u, (\nabla \mathbb{P}_0^t)^t~ \Phi ~ \nabla
\mathbb{P}_n^t\rangle =0.$$

Reciprocally, if relations (\ref{pip}) and (\ref{12b}) hold, then, $L^*[u] = 0$,
since
$$\langle L^*[u], \mathbb{P}_n^t \rangle = \langle u, L[\mathbb{P}_n^t]
\rangle =  \sum_{i=\max\{1,n-s\}}^{n+s}\langle u, \mathbb{P}_i^t\rangle
\Lambda_i^n =0.$$
Therefore, Lemma 3.4 in \cite{FPP1} gives
\begin{eqnarray*}
 0 &= &\langle L^*[u],
\mathbb{P}_{m}~ \mathbb{P}_n^t\rangle = \langle
u, L[\mathbb{P}_{m}~\mathbb{P}_{n}^t]\rangle \\
&=& \langle u, L[\mathbb{P}_{m}]~\mathbb{P}_{n}^t \rangle +
\langle u, \mathbb{P}_{m}~L[\mathbb{P}_{n}^t] \rangle + 2 \langle
u, (\nabla \mathbb{P}_{m}^t)^t ~ \Phi ~ \nabla
\mathbb{P}_n^t\rangle.
\end{eqnarray*}

If $0\le m< n-s$, $\langle u,
L[\mathbb{P}_{m}]~\mathbb{P}_{n}^t \rangle = 0$, since
$\deg(L[\mathbb{P}_m]) \le m+s$. On
the other hand, from (\ref{pip}), we deduce $\langle u,
\mathbb{P}_{m}~L[\mathbb{P}_{n}^t]\rangle =0$, and therefore, we conclude
$$ \langle u,
(\nabla \mathbb{P}_{m}^t)^t \Phi\, \nabla \mathbb{P}_n^t\rangle=0.$$
Finally, Theorem \ref{teor3} provides the desired
result.
\end{proof}

Again, three term recurrence relations (\ref{rr3t1}) and (\ref{rr3t2}), allow us to
express $L[\mathbb{P}_n^t]$ in terms of the vector polynomials
$\mathbb{P}_{n+1}^t$ and $\mathbb{P}_{n}^t$:

\begin{cor} If a
quasi--definite moment functional $u$
is semiclassical, then there exist polynomial matrices $N_1^n$, and $N_2^n$, satisfying
\begin{equation}
L[\mathbb{P}_n^t]=
\mathbb{P}_{n+1}^t~ N_{1}^n+ \mathbb{P}_{n}^t~ N_{2}^n
\end{equation}
with $\deg(N_1^{n})\leq s-1$ and $\deg(N_2^{n})\leq s$.
\end{cor}

\section{Examples}

\noindent{\bf Example 1: Appell--type orthogonal polynomials}

The so--called {\it Appell polynomials} (\cite{AK}) are orthogonal
polynomials in $d$ variables associated with the weight function
$$
\omega_{\alpha}(\textrm{x}) =
\textrm{x}^{\alpha}\,(1-|\textrm{x}|)^{\beta},$$ where $\alpha =
(\alpha_1,\dots,\alpha_d)\in \mathbb{R}^d$, with $\alpha_i> -1$,
$1\le i \le d$, and $\beta\in \mathbb{R}$ such that $\beta
>-1$,
on the $d$--simplex,
$$T_d=\{\textrm{x}=(x_1,\dots,x_d):\quad x_1\geq 0,\dots,x_d\geq
0, \quad  1-|\textrm{x}|\geq 0\}.$$

The Appell moment functional $u$ is defined as follows
$$\langle u, f\rangle = \int_{T_d}\, f(\textrm{x}) \,
\omega_{\alpha}(\textrm{x})\, d\textrm{x}.$$ Using our definitions
(\cite{AFPP1},\cite{FPP3}), $u$ is a classical moment functional
(i.e., it is semiclassical with $s=0$), since it satisfies the
matrix Pearson--type equation (\ref{tpm}), where the matrices
$\Phi$ and $\Psi$ are given by
$$\Phi = \begin{pmatrix}
                x_1(x_1-1) & x_1 x_2  & \cdots & x_1 x_d\cr
                x_2 x_1    & x_2(x_2-1) & \cdots & x_2 x_d \cr
                \vdots     & \vdots     & \ddots & \vdots \cr
                x_d x_1    & x_d x_2    & \cdots & x_d(x_d-1)\cr
                \end{pmatrix},$$
$$    \Psi =\begin{pmatrix}
           (|\alpha|+d)x_1-(\alpha_1+1) \cr
           (|\alpha|+d)x_2-(\alpha_2+1) \cr
            \vdots\cr (|\alpha|+d)x_d-(\alpha_d+1)\cr
            \end{pmatrix}.
$$
Now, we introduce the {\it Appell--type polynomials} as the
orthogonal polynomials in $d$ variables with respect to the moment
functional
$$v = u + \lambda \,
\delta(\textrm{x}),$$ where $\lambda \ge0$ is a positive real
number, and $\delta(\textrm{x})$ is the usual Dirac distribution at
${\bf 0}\in \mathbb{R}^d$. The action of $v$ over polynomials is
defined as follows,
$$\langle v, f\rangle = \int_{T_d}\, f(\textrm{x}) \,
\omega_{\alpha}(\textrm{x})\, d\textrm{x} + \lambda \, f({\bf
0}).$$

The moment functional $v$ is semiclassical with $s=1$,
since $v$ satisfies the matrix Pearson--type equation (\ref{tpm})
\begin{equation}\label{A--t}
\hbox{\rm div~} (\hat\Phi ~v) = \hat\Psi^t ~v,
\end{equation}
 where
$$\hat\Phi = x_1\,\Phi, \qquad \hat\Psi = (x_1(x_1-1), x_1 x_2, \cdots, x_1 x_d)^t + x_1\Psi.$$

In fact, using that $x_1\,\delta(\textrm{x}) = 0$, we get
\begin{eqnarray*}
\hbox{\rm div~} (\hat\Phi ~v) &=& \hbox{\rm div~} (x_1\,
\Phi ~(u + \lambda \delta(\textrm{x})) =\\
&=& \hbox{\rm div~} (x_1\,
\Phi ~u ) + \hbox{\rm div~} (x_1\,
\Phi ~\lambda ~\delta(\textrm{x})) =\\
&=& (1,0, \cdots,0)\, \Phi ~u + x_1 \, \hbox{\rm div~} (\Phi ~u) =\\
&=&(x_1(x_1-1), x_1 x_2, \cdots, x_1 x_d) v + x_1 \, \Psi^t \, v.
\end{eqnarray*}

Observe that the matrix Pearson--type equation for $v$ is not
unique. The moment functional $v$ also satisfies (\ref{A--t}) with
$$\hat\Phi = x_i \Phi, \quad \hat\Psi = (x_i x_1, \cdots, x_i(x_i-1), \cdots, x_i x_d)^t + x_i\Psi,$$
for $1\le i \le d$, since $x_i\,\delta(\textrm{x}) = 0$.

\noindent{\bf Example 2: A multivariate analogue of the classical orthogonal polynomials}

Examples of two--variables analogues of the Jacobi polynomials are
studied in \cite{Koo} by T. Koornwinder. Using  similar tools, we
present an example of a semiclassical weight function with
unbounded support.

Let $\alpha_i$ be real numbers with $\alpha_i>-1$, for $1\le
i\le d$. Then, for $k_i\ge 0$, $1\le i\le d$, and $k_1\ge k_d,$ we
define the polynomials

$$P^{(\alpha_1, \ldots, \alpha_d)}_{k_1, \ldots, k_d}(\textrm{x})=  L_{k_1-k_d}^{(\alpha_1+ 2 k_d +1)} (x_1)
L_{k_{2}}^{(\alpha_{2})} (x_2)\ldots L_{k_{d-1}}^{(\alpha_{d-1})} (x_{d-1})\,x_1^{k_d} P_{k_d}^{(\alpha_d, 0)} (x_1^{-1}
x_d), $$
where $L_{k_i}^{(\alpha_i)}$ is a Laguerre polynomial in one variable,
and $P_{k_d}^{(\alpha_d, 0)}$ is a Jacobi polynomial.

The polynomials $P^{(\alpha_1, \ldots, \alpha_d)}_{k_1, \ldots,
k_d}(\textrm{x})$ are orthogonal with respect to the weight
function
$$w(\textrm{x})=x_1^{\alpha_1}\ldots x_{d-1}^{\alpha_{d-1}}\, e^{-(x_1+\ldots +x_{d-1})}\, (1-x_1^{-1} x_d)^{\alpha_d},$$ on the
region  $\{\textrm{x}=(x_1,\ldots, x_d) \, /\,-x_1 < x_d
<x_1,\,\,\, x_i
>0, \, i=1,  \ldots, d-1\}$.

Defining the matrices
$$\Phi=\left(%
\begin{array}{ccccc}
  x_1(x_1-x_d) & &   &  & \bigcirc \\
   & x_2 &  & &  \\
 &  &   \ddots &  & \\
   & &  & x_{d-1} & \\
\bigcirc &  &  &  & x_1^2(x_1-x_d)\\
\end{array}%
\right),$$  $$\Psi=\begin{pmatrix}-x_1^2+x_1\, x_d+ (\alpha_1 +2)x_1
+ (\alpha_d-\alpha_1-1)x_d \cr \alpha_2-x_2 \cr \vdots \cr
\alpha_{d-1}-x_{d-1} \cr-(\alpha_d+1)x_1^2\end{pmatrix},$$ we can
prove that the weight function $w(\textrm{x})$ satisfies (\ref{tpm})
and so, the polynomials $P^{(\alpha_1, \ldots, \alpha_d)}_{k_1,
\ldots, k_d}(\textrm{x})$ are semiclassical.

Another examples for semiclassical orthogonal polynomials in two
variables appear in \cite{AFPP1}.

\end{document}